\documentclass[12pt]{article}
\usepackage{amsfonts}
\usepackage{amsmath}
\usepackage{amsthm}
\usepackage{amssymb}
\usepackage{mathrsfs}
\usepackage{graphicx}

\newtheorem{thm*}{Theorem}
\renewcommand{\iint}{\int \!\!\!\! \int}

\long\def\symbolfootnote[#1]#2{\begingroup%
\def\thefootnote{\fnsymbol{footnote}}\footnote[#1]{#2}\endgroup}
\title{Existence of Spherical Initial Data with Unit Mass, Zero Energy, and Virial less than - 1/2 for the Relativistic Vlasov-Poisson Equation with Attractive Coupling }

\begin{document}

\maketitle

\begin{abstract}
\noindent
In a recent paper, Kiessling and Tahvildar-Zadeh proved that any classical solution of the relativistic Vlasov-Poisson equation with attractive coupling launched by spherically symmetric initial data with zero total energy and virial less than or equal to -1/2 will blow up in finite time.  They left open whether such data exist.  Subsequently, the question was raised whether any such data exist at all.  In fact, the simplest conceivable ansatz, a nearly uniform ball of material centered at the origin with momenta directed inward, must have virial strictly larger than -1/2!  In this brief note, we settle this issue by constructing two classes of such initial data.
\end{abstract}

\section{Introduction}

The relativistic Vlasov-Poisson (rVP) system is given by $$\textrm{rVP}^{\pm}:\;\left\{ \begin{array}{r}\left(\partial_t + \frac{p}{\sqrt{1+\lvert p \rvert^2}} \cdot \nabla_q \pm \nabla_q\varphi_t(q) \cdot \nabla_p  \right)f_t(p,q)=0\\ \\ \triangle_q\varphi_t(q) = 4\pi \int f_t(p,q)\; d^3p \\ \\\varphi_t(q) \asymp -\lvert q \rvert^{-1} \textrm{ as } \lvert q \rvert \to \infty ; \end{array}\right.$$ $\textrm{rVP}^{+}$ models a system with repulsive interaction while $\textrm{rVP}^{-}$ models a system with attractive interaction.  One of the earliest papers to appear on the subject is \cite{GS85}  wherein Glassey and Schaeffer show that global classical solutions to $\textrm{rVP}^{\pm}$ will exist for initial data that are spherically symmetric, compactly supported in momentum space, vanish on characteristics with vanishing angular momentum, and have $\mathfrak{L}^{\infty}$-norm below a critical constant $\mathcal{C}_{\infty}^{\pm}$ with $\mathcal{C}_{\infty}^{+} = \infty$ and $\mathcal{C}_{\infty}^{-} < \infty.$  Recently, Kiessling and Tahvildar-Zadeh \cite{KTZ08} have extended the theorem of Glassey and Schaeffer for rVP$^-$ by proving global existence of classical solutions for the same type of initial data except that they only need to have $\mathfrak{L}^{\beta}$-norm below a critical constant $\mathcal{C}_{\beta}^{-}$ with  $\mathcal{C}_{\beta}^{-} < \infty$ and identically zero for $\beta < 3/2.$  The sharp value of this constant was subsequently determined for all values of $\beta$ in \cite{You10}.

Glassey and Schaeffer also investigated what may happen when rVP$^-$ is launched by initial data with $\|f\|_{\infty} > \mathcal{C}_{\infty}^{-}$.  They proved that negative energy data lead to ``blow-up" in finite time.   In the same paper extending the global existence result, Kiessling and Tahvildar-Zadeh  proved that any spherically symmetric classical solution of rVP$^-$ launched by initial data $f_0 \in \mathfrak{P}_3\cap\mathfrak{C}^1$ with \emph{zero total energy} and total virial less than or equal to -1/2 will blow up in finite time\symbolfootnote[2]{By $\mathfrak{P}_3\cap\mathfrak{C}^1$ we mean the set of probability measures on $\mathbb{R}^6$ whose first three moments are finite, that are absolutely continuous with respect to Lebesgue measure, and whose Radon-Nikodym derivative is $\mathfrak{C}^1$.}.  Specifically, with the total energy of $f_0$ given by
\begin{equation}
\mathcal{E}(f_0)=\int\!\!\! \sqrt{1+|p|^2}f(p,q)d\mu -\frac{1}{2} \iint \frac{f(p',q')f(p,q)}{\lvert q - q' \rvert} \;d\mu d\mu',\label{energy}
\end{equation}
and the virial of $f_0$ given by
\begin{eqnarray}
\mathcal{V}(f_0) &=& \int (q \cdot p)f_0(p,q)d\mu,\label{virial}
\end{eqnarray}
where $\int \cdots d\mu \equiv \int\!\!\int \cdots d^3pd^3q$, Theorem 6.1 of \cite{KTZ08} states:
\begin{thm*}
Let $t \to f_t \in \mathfrak{P}_3\cap\mathfrak{C}^1$ be a spherically symmetric classical solution of rVP$^-$ over some interval $[0, T)$ launched by initial data $f_0$ with $$ \| f_0 \|_{3/2} > \frac{3}{8}\left(\frac{15}{16}\right)^{1/3}. $$  If $\mathcal{E}(f_0) = 0$ and $\mathcal{V}(f_0) \le -1/2$, then $T < \infty$
\end{thm*}
The proof is by contradiction.  The zero-energy condition ensures that the virial is strictly decreasing.  Assuming that an initial datum satisfying the requirements above launches a global solution, one can show (using a number of technical estimates) that a strictly positive quantity is bounded above by a sum of positive terms and a time integral of the virial.  The authors then show that after a sufficiently long time this strictly positive quantity is bounded above by something negative.  Hence, such an initial datum cannot lead to global existence in time.

Now, the authors of \cite{KTZ08} did not give any specific examples of favorable initial data, and recently the question was raised whether such data actually exist \cite{Cal10}.  Indeed, nearly uniform balls of material centered at the origin (the simplest ansatz one could conceive) cannot be made to work.  Surprisingly, the existence of these initial data is an interesting, non-trivial question!  The purpose of this brief note is to settle this issue.

We first define a generic separation of variables ansatz and derive all basic formulae needed to carry out our subsequent investigations.  We then show that the simplest ansatz (a nearly uniform ball of material centered at the origin) cannot have the desired properties. Next, we employ a disjoint core-halo structure (a central ball of material with a separate thin ring of material further out) to construct a class of initial data satisfying the requirements in \cite{KTZ08}.  In addition, we show that zero-energy initial data can have arbitrarily negative virial.  The success of the core-halo scheme raises the question of whether such a disjoint structure is necessary for our requirements.  We answer this in the negative by exhibiting a second class of initial data with zero total energy and total virial less than or equal to -1/2 that have the nice properties of being monotonically decreasing and ``singly supported."

We should mention that the details of the blow-up scenario are not addressed in either \cite{GS85} or \cite{KTZ08}.  However, in \cite{LMR08} Lemou, M\'ehats, and Rapha\"el proved that ``blow-up" for systems launched by initial data with negative total energy approaches a self-similar collapse.

Lastly, we also mention other recent results on $\textrm{rVP}^-$. Had\v zi\'c and Rein \cite{HR07} showed the non-linear stability of a wide family of steady-state solutions of $\textrm{rVP}^-$ against a class of perturbations utilizing energy-Casimir functionals.  Shortly thereafter, Lemou, M\'ehats, and Rapha\"el \cite{LMR08,LMR09} investigated non-linear stability versus the formation of singularities in $\textrm{rVP}^-$ through concentration compactness techniques.  While these results are interesting in their own right, they will play no role in the following considerations.

\section{The Separation of Variables Ansatz}

We first establish the various formulae we shall need for our particular class of functions. Consider a generic separation-of-variables ansatz:
\begin{eqnarray}
f(p,q) &=& \mathcal{C}\eta(|q|)\Phi(|p|)\mathscr{L}(\cos(\theta_{p,q})),\label{def}
\end{eqnarray}
where $\theta_{p,q}$ is the angle between $q$ and $p$ (both considered as vectors in $\mathbb{R}^3$) and $\eta, \Phi,$ and $\mathscr{L}$ are non-negative.  The coefficient $\mathcal{C}$ will be chosen to normalize the total mass to $1$.

The total mass is given by
\begin{align*}
\int f(p,q)d\mu = &\mathcal{C}\!\! \int_0^{\infty}\!\!\!\!\!\eta(|q|)|q|^2d|q| \!\! \int_0^{\infty}\!\!\!\!\Phi(|p|)|p|^2d|p|
\iint_{\mathbb{S}^2 \times \mathbb{S}^2}\!\!\!\!\!\!\!\!\mathscr{L}(\cos(\theta_{p,q}))d\Omega_pd\Omega_q\\
=&\; 8\pi^2\mathcal{C} \int_0^{\infty}\eta(|q|)|q|^2d|q| \int_0^{\infty}\Phi(|p|)|p|^2d|p|\int_{-1}^1\mathscr{L}(x)dx.
\end{align*}
Using the notation $$ \| g \|_1 \equiv \int_0^{\infty} |g(r)| dr,  $$ and recalling that all functions under consideration are non-negative we see that taking
 $$\mathcal{C}^{-1} =  8\pi^2 \| \eta |q|^2 \|_1 \| \Phi |p|^2 \|_1 \int_{-1}^1\mathscr{L}(x)dx$$
will ensure that $f$ has total mass 1.  Similarly, the $\mathfrak{L}^{3/2}$-norm is given by
\begin{eqnarray*}
\|f\|_{3/2} &=& \frac{\left(\int_0^{\infty}\eta(|q|)^{3/2}|q|^2d|q| \int_0^{\infty}\Phi(|p|)^{3/2}|p|^2d|p|\int_{-1}^1\mathscr{L}(x)^{3/2}dx\right)^{2/3}}{2\pi^{2/3} \| \eta |q|^2 \|_1 \| \Phi |p|^2 \|_1 \int_{-1}^1\mathscr{L}(x)dx}.
\end{eqnarray*}

Recalling the definition of the virial functional \eqref{virial} from the introduction, for our separation of variables gives
\begin{eqnarray}
\mathcal{V}(f) &=& \mathcal{C}\!\!\int_0^{\infty}\!\!\!\!\!\eta(|q|)|q|^3 d|q|\int_0^{\infty}\!\!\!\!\!\Phi(|p|)|p|^3 d|p| \iint_{\mathbb{S}^2 \times \mathbb{S}^2} \!\!\!\!\!\!\!\!\cos(\theta_{p,q})\mathscr{L}(\cos(\theta_{p,q})) d\Omega_p d\Omega_q \nonumber\\
&=& 8\pi^2\mathcal{C} \| \eta |q|^3 \|_1 \| \Phi |p|^3 \|_1 \int_{-1}^1x\mathscr{L}(x)dx\nonumber\\
&=& \frac{\| \eta |q|^3 \|_1}{\| \eta |q|^2 \|_1}\frac{\| \Phi |p|^3 \|_1}{\| \Phi |p|^2 \|_1}\frac{\int_{-1}^1x\mathscr{L}(x)dx}{\int_{-1}^1\mathscr{L}(x)dx}.
\end{eqnarray}
Of course, in order to have a negative virial we must require $$\int_{-1}^1x\mathscr{L}(x)dx < 0.$$

For the energy (as given by \eqref{energy}), we compute the ``kinetic" and potential energy contributions separately (the portion of the energy we label as ``kinetic" also contains the rest mass energy).  For the kinetic energy, we have
\begin{eqnarray}
KE(f)\!\!\!&\equiv&\!\!\! \int  \sqrt{1+|p|^2} f(p,q) d\mu\nonumber\\
 &=&\!\!\!\mathcal{C}\!\!\! \int_0^{\infty}\!\!\!\!\Phi(|p|)\sqrt{1+|p|^2}|p|^2d|p|
\int_0^{\infty}\!\!\!\!\!\! \eta(|q|)|q|^2 d|q| \iint_{\mathbb{S}^2 \times \mathbb{S}^2} \!\!\!\!\!\!\!\!\!\!\!\!\mathscr{L}(\cos(\theta_{p,q})) d\Omega_p d\Omega_q \nonumber\\
&=&\!\!\! \frac{\|\Phi \sqrt{1+|p|^2}|p|^2\|_1}{\|\Phi |p|^2 \|_1}.
\end{eqnarray}

We begin computing the potential energy by first computing the spatial distribution $\rho$ associated to $f$:
\begin{eqnarray*}
\rho(q) &\equiv& \int f(p,q) \;d^3p\\
&=& \mathcal{C}\eta(|q|)\int_0^{\infty}\Phi(|p|)|p|^2d|p|\int_{\mathbb{S}^2}\mathscr{L}(\cos(\theta_{p,q})) d\Omega_p\\
&=& \frac{\eta(|q|)}{4\pi\| \eta|q|^2 \|_1}
\end{eqnarray*}
In terms of this function, the potential energy contribution is now given by $$PE(f) \equiv -\frac{1}{2}\iint \frac{\rho(q)\rho(q')}{|q-q'|}\; d^3q' \; d^3q.$$
We calculate this integral directly, and after a little work find
\begin{eqnarray}
PE(f) &=& -\frac{1}{\| \eta|q|^2 \|_1^2}\int_0^{\infty} \eta(|q|)|q| \left(\int_0^{|q|}\eta(|q'|)|q'|^2d|q'| \right)d|q|.
\end{eqnarray}

\section{Nearly Uniform Balls Fail (Barely)}

At this point, we can show that the simplest ansatz cannot work.  Relaxing the differentiability requirement for a moment, we make the specific choices
\begin{eqnarray*}
\eta(|q|) &=& \chi_{[0,R]}(|q|), \\
\Phi(|p|) &=& \chi_{[0,P]}(|p|), \\
\mathscr{L}(x) &=& \chi_{[-1,a]}(x)
\end{eqnarray*}
where $R>0, P>0, -1<a < 1,$ and $\chi_{I}$ is the characteristic function of the interval $I$. We have that
\begin{eqnarray*}
KE(f) &=& \frac{3}{8}\left(\frac{\sqrt{1+P^2}}{P^2} + 2\sqrt{1+P^2} - \frac{\ln(P+\sqrt{1+P^2})}{P^3}  \right),\\
PE(f) &=& -\frac{3}{5R}.
\end{eqnarray*}

The zero-energy requirement allows us to solve for $R$ in terms of $P$.  Calculating the virial and plugging in the formula for $R = R(P)$ gives $\mathcal{V}(f)$ in terms of the parameters $P$ and $a$.  Simple asymptotics shows that for this ansatz $$\mathcal{V}(f) > -\frac{9}{20}$$ for any choice of parameters.  Smoothing out the boundaries of these step functions (taking care to keep the total energy zero) shows that this ansatz is untenable.  Of course, this type of initial condition may well lead to collapse after finite time, but this cannot be verified by our virial condition.

The failure of nearly uniform balls to satisfy our criteria is a little surprising at first, but upon reflection we can understand the difficulty.  Recall that the ``kinetic" energy term contains the rest mass contribution and so is at least one (in our rationalized units).  Heuristically, this forces the radius of the ball of mass to be small in order to have the potential energy negative enough to balance the ``kinetic" energy contribution.  This has severe consequences for the virial - which is made smaller in magnitude by concentrating all the matter near the origin.

\section{Success: A Disjoint Core-Halo Ansatz}

The considerations at the end of the previous section lead us to try a disjoint core-halo scheme:
\begin{eqnarray}
\eta(|q|) &=& \chi_{[0,R_1]}(|q|) + \alpha\chi_{[R_2,R_3]}(|q|), \label{sp1}\\
\Phi(|p|) &=& \chi_{[0,P]}(|p|), \label{sp2}\\
\mathscr{L}(x) &=& \chi_{[-1,a]}(x) \label{sp3}
\end{eqnarray}
where $0< R_1 \le R_2 \le R_3, 0 < \alpha, 0 < P, $ and $-1 < a \le 1$.  Again, we have relaxed the differentiability requirement to make the computations tractable.  The idea behind this ansatz is that putting a tiny halo of material far from the center should increase the magnitude of the virial while keeping the potential energy negative enough to maintain zero total energy for the system.

We first consider the specific choices $$R_1 = \frac{1}{5},\quad R_2 = 1,\quad R_3 = 2,\quad \textrm{ and } P = 1  .$$  The zero energy condition forces (thanks to Maple) $$\alpha = \frac{35\ln(1+\sqrt{2})+30 - 105\sqrt{2} + 2\sqrt{6480\sqrt{2}-1655-2160\ln(1+\sqrt{2})}}{125\left(735\sqrt{2} - 188 - 245 \ln(1+\sqrt{2})\right)} $$
(a small but positive number).  Again by Maple, choosing $-1 < a \le -4/5$ gives $\mathcal{V}(f) < -1/2.$  Since our ansatz is compactly supported, we have the requisite number of moments.  To complete the argument, we note that we can smooth out the various step functions in such a way that the resulting integrals are as close to the values obtained above as we like.  Since $\alpha$ was chosen to force the zero-energy condition and the necessary value was strictly positive, it follows that by smoothing the test functions appropriately we can still choose an $\alpha > 0$ to keep the total energy zero.  Using Maple to check the $\mathfrak{L}^{3/2}$-norm for the range of $a$ listed above shows that our data are well above the critical norm $\frac{3}{8}\left(\frac{15}{16}\right)^{1/3}.$ Thus, we have constructed an initial datum satisfying all the requirements of Theorem 6.1 in \cite{KTZ08}.

Considering the lower bound on the virial for nearly uniform balls, the reader may wonder whether a similar phenomenon holds for the disjoint core-halo structure.  We give an argument showing that the virial can be made arbitrarily negative in this case.  We make the following choices in (\ref{sp1}):
$$R_1 = P^{-2},\quad R_2 = P,\quad \textrm{ and } R_3 = P^2.$$ Now, the zero energy condition allows us to solve for $\alpha = \alpha(P)$.  Using Maple for convenience, one sees that as a function of $P$, $\alpha$ is positive for sufficiently large $P$ and is proportional to $P^{-23/2}$ as $P \to \infty$.  Plugging these choices into the formula for the virial and looking at the asymptotics for large $P$ shows that the virial is proportional to $-(1-a)P^3$. Smoothing these functions out as above shows that the virial is unbounded below.  As before, the $\mathfrak{L}^{3/2}$-norm for these functions remains well above the critical norm for a wide range of values of the parameter $a$.

\section{A Monotonic Core-Halo Ansatz}

Having established the existence of initial data satisfying the theorem of Kiessling and Tahvildar-Zadeh, we are technically done.  However, we can ask whether the separation of the halo from the core is an absolute necessity.  The answer is no.  Even though nearly uniform balls (after regularization) were shown not to work, we can find zero-energy initial data with virial less than $-1/2$ that are monotonically decreasing and supported in a ball centered at the origin.

As in the previous section, we describe initial data that are not smooth with the understanding that we can appropriately regularize our functions without significantly altering any of the required properties.  We now consider the ansatz
\begin{eqnarray}
\eta(|q|) &=& \chi_{[0,R_1]}(|q|) + \left(\frac{R_1}{|q|}\right)^n\chi_{[R_1,R_2]}(|q|) + \left(\frac{R_1}{R_2}\right)^n\chi_{[R_2,R_3]}(|q|), \label{sp4}\\
\Phi(|p|) &=& \chi_{[0,P]}(|p|), \label{sp5}\\
\mathscr{L}(x) &=& \chi_{[-1,a]}(x). \label{sp6}
\end{eqnarray}
Note that we have kept the same choices for the momentum and angular portions of $f$ but substituted a spatial portion that is monotonically decreasing inside the ball of radius $R_3$ centered at the origin and zero outside this ball.  A representative diagram for $\eta$ is given in Figure \ref{fig3} below.  Note that it retains the essential ingredient from the previous section.  Namely, there is a dense ball of material centered about the origin (to sufficiently lower the potential energy) surrounded by a thin atmosphere of material further out (to sufficiently increase the magnitude of the virial).

\begin{figure}[ht]\centering
  \includegraphics[width=75mm]{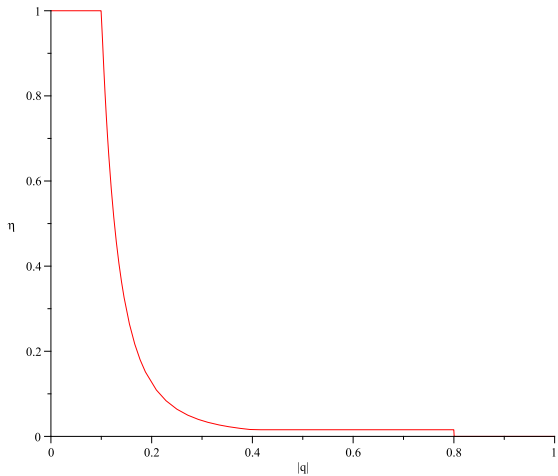}\\
  \caption{Representative monotonically decreasing $\eta$}
  \label{fig3}
\end{figure}

Even with this relatively simple ansatz, the calculations are unwieldy.  Thanks again to Maple, we merely report a choice of parameters that does the job.  Namely, choosing $$R_1 = \frac{1}{100},\quad R_2 = \frac{1}{11},\quad R_3 = \frac{1}{10},\quad n = 3  $$ sets a corresponding $P$ (about $19.69$) to force the zero-energy condition.  Finally, any choice of $a$ less than roughly $-9/10$ gives us a virial less than $-1/2$.  A check of the $\mathfrak{L}^{3/2}$-norm shows that it is well above the critical value.  We remark that our choice of parameters was totally arbitrary, and so we expect that an entire range of such parameters will do the job.  However, here we are content merely with the existence of one such solution!

\section{Acknowledgements}
The author wishes to thank Simone Calogero for raising the issue.  The author is also grateful to Shadi Tahvildar-Zadeh and Michael Kiessling for many enlightening conversations!  The author also gratefully acknowledges funding through NSF grant DMS 08-07705 awarded to Michael Kiessling.

\end{document}